\documentclass[11pt]{article}
\usepackage{epsfig}
\title{Dynamics of a Nonlocal Kuramoto-Sivashinsky Equation
\footnote{This work was supported  by the National Science
Foundation Grant DMS-9704345.}  }
\author{  }
\date{August 5, 1997}
\begin{document}

\maketitle

\begin{center}
 Jinqiao Duan and Vincent J. Ervin \\
Department of Mathematical Sciences\\
Clemson University\\
Clemson, SC 29634-1907\\
USA \\
E-mail: duan@math.clemson.edu  \\
\hspace*{0.52in} ervin@math.clemson.edu
\end{center}

\newpage

{\bf Running Head}: Nonlocal Kuramoto-Sivashinsky Equation

{\em Note: This paper is in Journal of Didderential Euqations,
	Vol. 143 (1998), 243-266.}
 
\medskip
\medskip
 
{\bf Address correspondence to}:
\begin{center}
        Jinqiao Duan
 
Department of Mathematical Sciences
 
Clemson University
 
Clemson, SC 29634-1907

USA
 
E-mail: duan@math.clemson.edu

Phone: (864) 656-2730
 
Fax: (864) 656-5230
\end{center}
 
\newpage

\newtheorem{theorem}{Theorem}[section]
\newtheorem{conjecture}{Conjecture}[section]
\newtheorem{corollary}{Corollary}[section]
\newtheorem{lemma}{Lemma}[section]
\newtheorem{definition}{Definition}[section]
\newtheorem{proposition}{Proposition}[section]
\newtheorem{problem}{Problem}[section]

\newcommand{\e}{\epsilon}
\renewcommand{\a}{\alpha}
\renewcommand{\b}{\beta}
\renewcommand{\aa}{\mbox{$\cal A$}}
\renewcommand{\L}{\mbox{$\cal L$}}

\newcommand{\integer}{\mbox{$\mathrm{Z\!\!Z}$}}
\newcommand{\real}{\mbox{$\mathrm{I\!R}$}}
\newcommand{\Natural}{\mbox{$\mathrm{I\!N}$}}
\newcommand{\bfx}{\mbox{$\mathbf{x}$}}
\newcommand{\bfu}{\mbox{$\mathbf{u}$}}
\newcommand{\bfv}{\mbox{$\mathbf{v}$}}
\newcommand{\bfn}{\mbox{$\mathbf{n}$}}
\newcommand{\bfb}{\mbox{$\mathbf{b}$}}
\newcommand{\bfc}{\mbox{$\mathbf{c}$}}
\newcommand{\bfr}{\mbox{$\mathbf{r}$}}
\newcommand{\bfxi}{\mbox{$\mathbf{\xi}$}}
\newcommand{\bfzeta}{\mbox{$\mathbf{\zeta}$}}
\newcommand{\mcR}{\mbox{$\mathcal{R}$}}

\newcommand{\al}{\alpha}
\newcommand{\be}{\beta}
\newcommand{\om}{\omega}
\newcommand{\La}{\Lambda}
\newcommand{\la}{\lambda}
\newcommand{\de}{\delta}
\newcommand{\beq}{\begin{equation}}
\newcommand{\eeq}{\end{equation}}
\newcommand{\beqo}{\begin{equation*}}
\newcommand{\eeqo}{\end{equation*}}
\newcommand{\p}{\partial}
\def\qed{\hbox{\vrule width 6pt height 6pt depth 0pt}}
\renewcommand{\theequation}{\thesection.\arabic{equation}}
\renewcommand{\thefigure}{\thesection.\arabic{figure}}
\renewcommand{\thetable}{\thesection.\arabic{table}}
\renewcommand{\thetheorem}{\thesection.\arabic{theorem}}
\catcode`\@=11
\@addtoreset{equation}{section}

\end{document}